\pdfoutput=1
\documentclass[11pt]{amsart}
\usepackage{geometry}                
\usepackage[parfill]{parskip}   
\usepackage{graphicx,color,amsthm}
\usepackage{amssymb}
\usepackage{epstopdf}
\usepackage{pdfsync}
\usepackage{caption}
\usepackage{subcaption}
\usepackage{mathtools}
\usepackage{fullpage}
\theoremstyle{plain} 
\newtheorem{thm}{Theorem}

\theoremstyle{definition} 
\newtheorem{defn}{Definition}
\newtheorem{ex}{Example}

\title{Computing Ribbon Obstructions for Colored Knots}
\author{Patricia Cahn and Alexandra Kjuchukova}
\thanks{This work was partially supported by the Simons Foundation/SFARI (Grant Number 523862, P. Cahn) and by NSF-DMS grants 1821212 to P. Cahn and 1821257 to A. Kjuchukova}
\begin{document}

\begin{abstract} 
Kjuchukova's $\Xi_p$ invariant gives a ribbon obstruction for Fox $p$-colored knots.  The invariant is derived from dihedral branched covers of 4-manifolds, and is needed to calculate the signatures of these covers, when singularities on the branching sets are present. In this note, we give an algorithm for evaluating $\Xi_p$ from a colored knot diagram, and compute a couple of examples.  \end{abstract}
\maketitle
MSC classes: 57M12, 57M25, 57Q60.

\section{Introduction}

A knot $\alpha\subset S^3$ is {\it slice} if it bounds an embedded disk in $B^4$; {\it ribbon} if it bounds such a disk  which has only local minima and saddle points with respect to the radial height function on $B^4$; and {\it homotopy ribbon}  if $K$ bounds a disk $D\subset B^4$ such that  the inclusion-induced map $\pi_1(S^3-\alpha)\twoheadrightarrow \pi_1(B^4-D)$ is surjective.  Every ribbon knot is homotopy ribbon, and every homotopy ribbon knot is slice.  The notions of slice and homotopy ribbon make sense in both the smooth and topological categories; the associated disk in $B^4$ is assumed to be smoothly or locally flatly embedded, respectively.  Ribbonness makes sense only in the smooth category.  Fox's Slice Ribbon Conjecture asserts that every  smoothly slice knot is ribbon~\cite{fox1962problems}.  The analogous assertion in the topological category would be that every topologically slice knot is homotopy ribbon.

Now suppose $\alpha$ is a Fox $p$-colored knot in $S^3=\partial X$, where $X$ is an oriented topological 4-manifold.  Let $\rho:\pi_1(S^3-\alpha)\twoheadrightarrow D_p$ denote the $p$-coloring of $\alpha$.   Kjuchukova's invariant $\Xi_p(\alpha,\rho)$ is defined for any colored knot $\alpha$ whose $p$-coloring extends over some locally flat, embedded, oriented surface $F\subset X$ with $\partial F=\alpha$ \cite{ geske2018signatures, kjuchukova2018dihedral}. When $\rho$ extends over a homotopy-ribbon disk for $\alpha$, the value  $\Xi_p(\alpha, \rho)$ falls in a bounded range.  In particular, when the $p$-fold dihedral branched cover of $S^3$ along $\alpha$ is a rational homology sphere, $|\Xi_p(\alpha,\rho)|\leq \dfrac{p-1}{2}$; in the general case an additional term appears in this inequality~\cite{cahnkjuchukova2017singbranchedcovers, geske2018signatures}.  In particular, when $p=3$, $|\Xi_3(\alpha,\rho)|=1$.  Hence, $\Xi_p$ provides a means for testing potential counterexamples to the Slice Ribbon Conjecture. This is one motivation to develop tools for the evaluation of this invariant, as the formula for $\Xi_p$ derived in~\cite{ geske2018signatures, kjuchukova2018dihedral} is not combinatorial or diagrammatic in nature. In addition to the purely knot-theoretic interest of $\Xi_p$, this procedure for evaluating the invariant also allows us to compute signatures of dihedral covers of four-manifolds with singular branching sets. Indeed, the invariant $\Xi_p(\alpha,\rho)$ was originally defined as the contribution to the signature of a dihedral cover $f: Y^4\to X^4$ whose branching set is embedded in $X$ with a singularity whose link is $\alpha$.

In this note, we lay out an algorithm for computing the ribbon obstruction $\Xi_p(\alpha,\rho)$ from a diagram of $\alpha$.  We then evaluate the value of $\Xi_p(\alpha,\rho)$ in several examples.  We focus on the case where $p=3$ for ease of exposition, but the procedure presented generalizes to all odd $p$.  That is, the computational algorithm given applies to all Fox $p$-colored knots for which the corresponding branched cover of $S^3$ along $\alpha$ is a rational homology sphere. 

Our main result, Theorem~\ref{procedurethm}, allows us to evaluate $\Xi_p(\alpha,\rho)$ from a colored diagram of $\alpha$ using a formula for $\Xi_p(\alpha,\rho)$ proved in~\cite{kjuchukova2018dihedral} (see Equation~\ref{eqXi}) as well as the algorithm developed in~\cite{cahnkjuchukova2016linking} for computing linking numbers of rationally null-homologous knots in dihedral covers of $S^3$. We apply this theorem in Section ~\ref{examples}.

$D_p$ denotes the dihedral group of order $2p$. In this paper, $p>1$ is an odd integer. 

\section{The invariant $\Xi_p$}

 $\Xi_p$ is an invariant of a knot $\alpha$ together with a choice of Fox $p$-coloring $\rho: \pi_1(S^3-\alpha)\twoheadrightarrow D_p$.  Rather than write $\Xi_p(\alpha,\rho)$, we write $\Xi_p(\alpha)$, as the choice of coloring is often understood or, up to equivalence of colorings, unique.

The invariant $\Xi_p(\alpha)$ arises in the following context.  Let $X$ be a closed, oriented 4-manifold.  Suppose that $B\subset X$ is a surface, embedded locally flatly away from one singular point whose link is $\alpha$. Given a surjection $\varphi: \pi_1(X-B)\twoheadrightarrow D_p$, we consider the induced {\it irregular dihedral $p$-fold cover} $Y$ of $X$ branched along $B$. This cover is characterized by the fundamental group of its unbranched counterpart, which is isomorphic to $\varphi^{-1}(\mathbb{Z}/2\mathbb{Z})$ for a choice of reflection subgroup in $D_p$.

 In the above setting, the invariant $\Xi_p(\alpha)$ should be viewed as the contribution to $\sigma(Y)$, the signature of the covering space, resulting from the presence of a singularity, $\alpha$, on the branching set.  For this reason, we refer to $\Xi_p$ as the signature defect associated to $\alpha$. Precisely, $\sigma(Y)$, $\sigma(X)$, and $\Xi_p$ are related as follows (see \cite[Theorem 1.4 (2)]{kjuchukova2018dihedral}):
 
$$\sigma (Y)=p\sigma (X)-\dfrac{p-1}{4}e(B)-\Xi_p(\alpha), $$
where $e(B)$ denotes the self-intersection number of the branching set.  The fact that $\Xi_p$ gives a ribbon obstruction when $Y$ is a manifold is proven in~\cite{cahnkjuchukova2017singbranchedcovers}. This obstruction is generalized to a larger class of knots, namely all colored knots which bound colored surfaces in some four-manifold with $S^3$ boundary, in~\cite{geske2018signatures}. Knots with this property are called  $p$-{\it admissible}.  Theorem~\ref{procedurethm} applies to all 3-admissible knots whose irregular 3-fold dihedral covers are rational homology spheres.

\begin{defn}
\label{chark}
Let $\alpha\subset S^3$ be a knot and $V$ a Seifert surface for $\alpha$ with Seifert form $A_V$. Let $L_V=A_V+A_V^{T}$ be the corresponding symmetrized form.  Let $\beta\subset V^\circ$ be an embedded curve representing a primitive class in $H_1(V; \mathbb{Z})$. If $L_V(\beta,\omega)=A_V(\beta, \omega) + A_V(\omega, \beta) \equiv 0\mod p$ for all embedded curves $\omega$ in $V$, we say that $\beta$ is a {\it mod~p characteristic knot} for $\alpha$.
\end{defn}

Characteristic knots are key for computing the signature defect $\Xi_p(\alpha)$. The existence of  a $p$-fold irregular dihedral cover of $S^3$ branched along $\alpha$ is equivalent to the existence of a {\it mod~p characteristic knot} for $\alpha$~\cite{CS1984linking}; the role of a characteristic knot is discussed further in Section~\ref{dih-con}.

\section{Overview of the Algorithm}\label{algoverview.sec}

\subsection{The signature defect arising from a singularity}
Our combinatorial procedure  for computing $\Xi_p(\alpha)$ relies on the formula given in  Theorem~1.3 of~\cite{kjuchukova2018dihedral}, which we now recall. Let $\alpha$ be a $p$-admissible knot,  $\rho: \pi_1(S^3 - \alpha)\to D_p$ a surjective homomorphism and $V$ a Seifert surface for $\alpha$.  Cappell and Shaneson showed~\cite[Proposition~1.1]{CS1984linking}, using the HNN presentation for $\pi_1(S^3 - \alpha)$, that the homomorphism $\rho$ can be described by linking curves in $S^3 - V$ with a characteristic knot  (see Definition~\ref{chark}) for $\alpha$ contained in the interior of $V$. Let $L_V$ denote the symmetrized Seifert form of $V$, $\beta$ the characteristic knot, $\zeta$ a primitive $p$-th root of unity, and $\sigma_{\zeta^i}$ the Tristram-Levine  $\sigma_{\zeta^i}$ signature~\cite{levine1969knot, tristram1969some}. We have, 
\begin{equation}\label{eqXi}
\Xi_p(\alpha) = \frac{p^2-1}{6p}L_V(\beta, \beta)  + \sum_{i=1}^{p-1} \sigma_{\zeta^i}(\beta) + \sigma(W(\alpha, \beta))
\end{equation}

 The first two terms in the above expression for $\Xi_p(\alpha)$ are easily calculated. The third term, $\sigma(W(\alpha, \beta))$, denotes the signature of a four-manifold $W(\alpha, \beta)$ constructed by Cappell and Shaneson in~\cite{CS1984linking}.   That is, $W(\alpha,\beta)$ is a cobordism between the $p$-fold irregular dihedral cover of $S^3$ branched along $\alpha$ and the $p$-fold cyclic cover of $S^3$ branched along $\beta$. We recall this construction in Section~\ref{octopus}. Computing the signature of the manifold  $W(\alpha, \beta)$ in terms of $\alpha$ is the main result of this paper, Theorem~\ref{procedurethm}, and this is equivalent to computing $\Xi_p(\alpha)$. The intersection matrix of $W(\alpha, \beta)$ can be expressed in terms of linking numbers of certain curves in the irregular $p$-fold dihedral cover $M_\alpha$ of $S^3$ branched along $\alpha$. We choose an orientation of $\beta$, and let $\beta_r$ and $\beta_l$ denote its right and left push-offs in $V$.  Let $V-\beta$ denote the surface with three boundary components $\alpha$, $\beta_r$, and $\beta_l$, obtained from $V$ by removing an annular neighborhood of $\beta$.  The curves whose linking numbers appear in the computation of $\sigma(W(\alpha, \beta))$ are lifts to $M_\alpha$ of a basis $\mathcal{B}_V:=\{\beta_r,\beta_l, \omega_1, ..., \omega_{2g-2}\}$  for $H_1(V-\beta; \mathbb{Z})$. The relevant linking numbers are computed using the algorithm given in~\cite{cahnkjuchukova2016linking}. 

We condense all this information in a labeled link diagram of $\alpha$, $\beta$ and the $\omega_j$, so that the signature defect can be computed algorithmically. The resulting algorithm is the content of Theorem~\ref{procedurethm}.

We set $p=3$ for the remainder of this section.

First, the arcs of $\alpha$ are labeled `1', `2', and `3', to indicate their coloring by the transpositions $(23)$, $(13)$, and $(12)$.  

\begin{figure}[htbp]\includegraphics[width=2in]{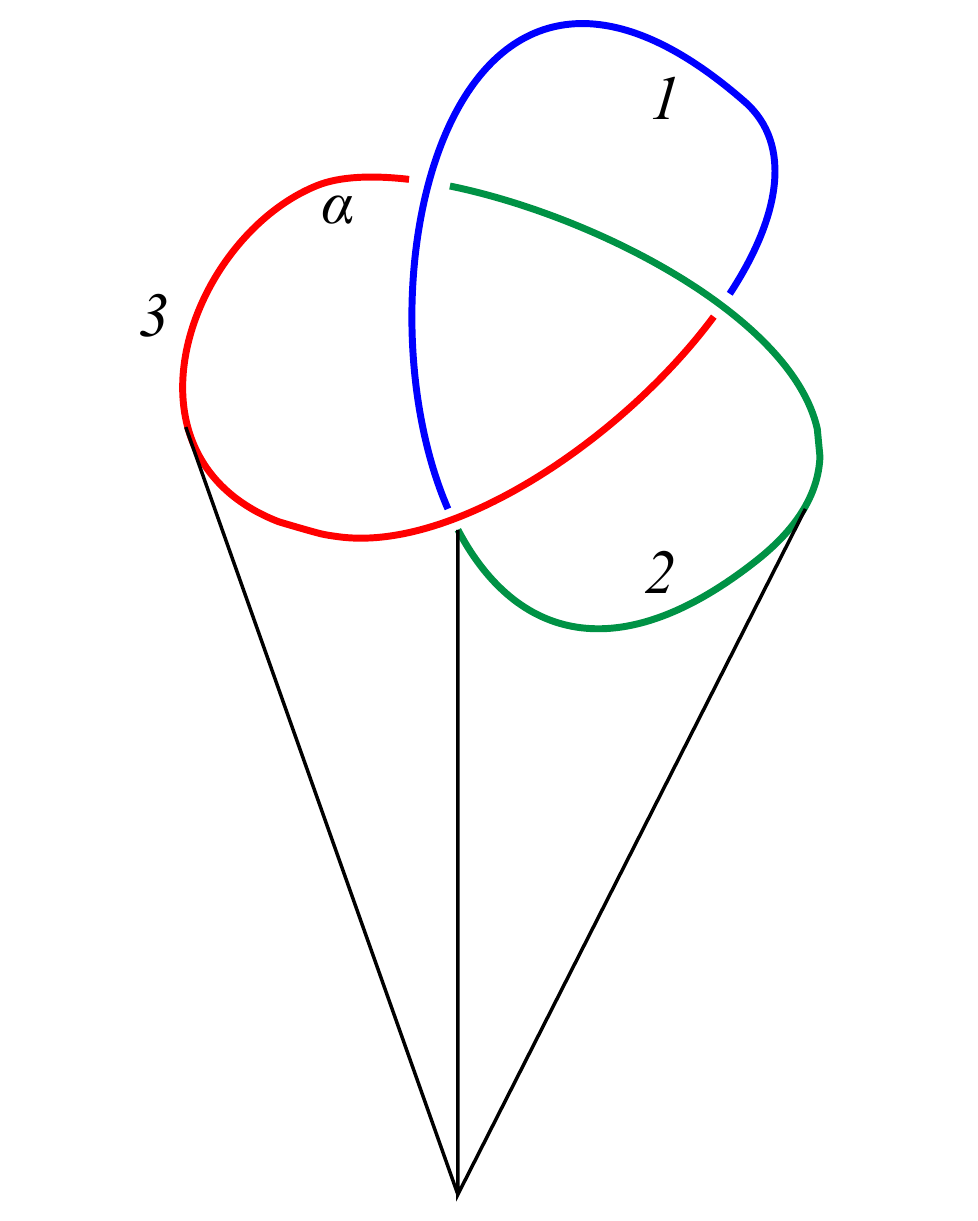}
	\caption{A cell structure on $S^3$ determined by the cone on $\alpha$.}
\end{figure}

\begin{figure}[htbp]\includegraphics[width=4in]{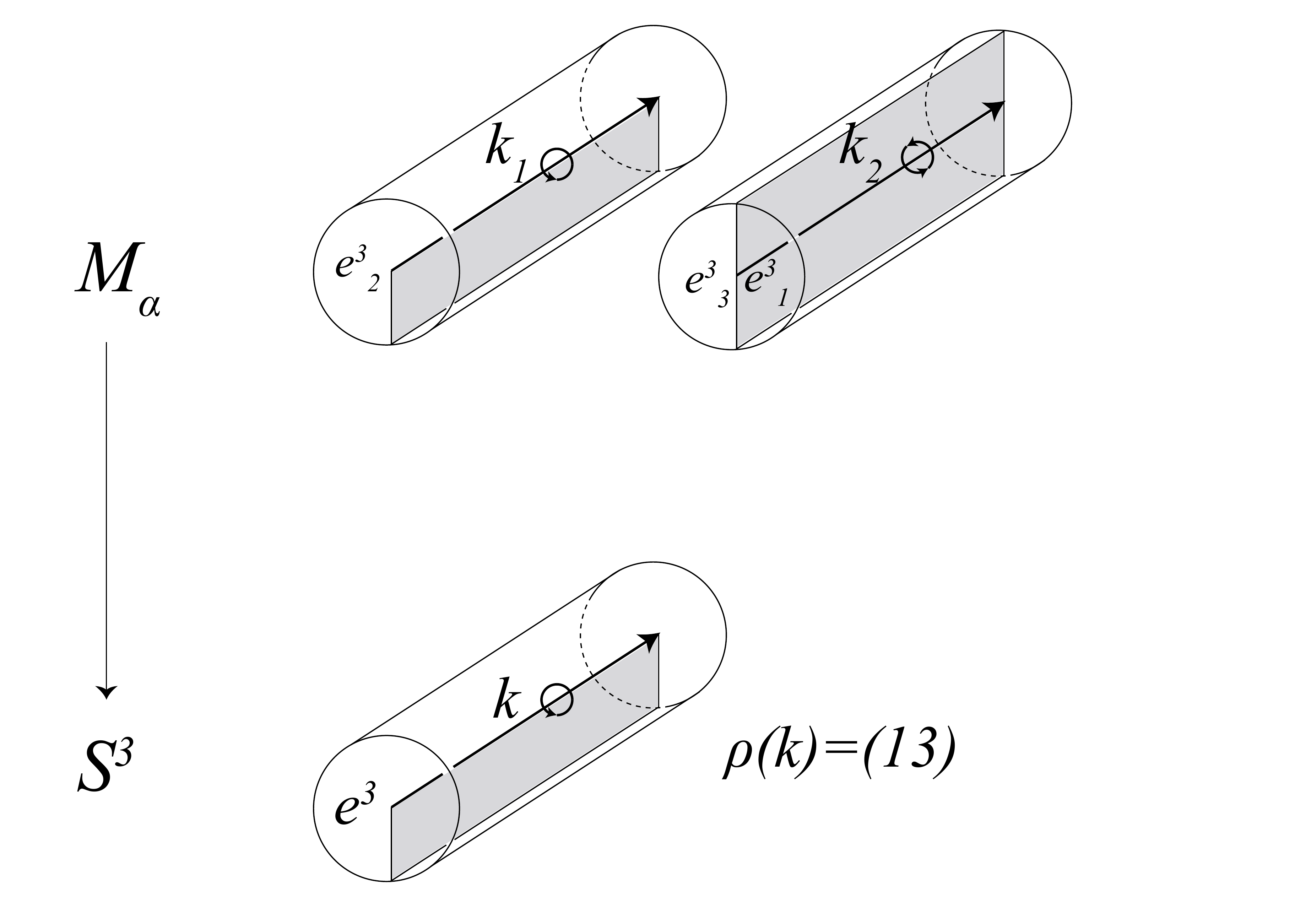}
\caption{Below, a tubular neighborhood of $k$, an arc in a diagram of $\alpha$. On top, the two lifts of this neighborhood, and the action of a meridian of $k$ on the three 3-cells, $e^3_1, e^3_2, e^3_3$,  of $M_\alpha$. The shaded rectangles represent 2-cells. Put together, these local pictures determine a cell structure on the 3-fold irregular dihedral cover of $S^3$, branched along $\alpha$.}\label{meridianliftprocedureintro.fig}
	
\end{figure}

The irregular dihedral cover $M_\alpha$ corresponding to $\rho$ is equipped with a cell structure determined by the cone on $\alpha$.  We review the key aspects of this cell structure here (see, e.g., \cite{cahnkjuchukova2016linking} for more details). First equip $S^3$ with a cell-structure that has one 3-cell $e^3$, the complement of the cone on $\alpha$.  The ``walls" of the cone on $\alpha$ are 2-cells, and so-on.  The cell structure on $M_\alpha$ is the lift of this cell-structure on $S^3$.  The 3-cell $e^3$ has three preimages in $M_\alpha$,  $e^3_1$, $e^3_2$, and $e^3_3$. These 3-cells are labeled such that the meridians of the arcs of $\alpha$ act on the subscripts according to the coloring $\rho:\pi_1(S^3-\alpha)\twoheadrightarrow D_3$.  See Figure ~\ref{meridianliftprocedureintro.fig}.

 We also choose a designated ``zeroth" arc of each component of the link diagram (later, all arcs in the link diagram will be numbered).   Let $\omega_i^j$ denote the lift of $\omega_i$ such that the lift of its zeroth arc lies in the $j^{th}$ 3-cell $e^3_j$, for $j=1,2,3$.  The lifts $\beta_r^j$ and $\beta_l^j$ are defined analogously. 

An {\it anchor path} for a curve $\omega\subset V-\beta$ is a properly embedded path $\delta$ in $V-\beta$ from a point $q$ on the zeroth arc of $\alpha$ to a point $r$ on the zeroth arc of $\omega$.  Suppose $\delta$ crosses under the arcs $a_1,\dots a_k$  in the diagram of $\alpha$, in that order, when traversing from $q$ to $r$.  The {\it monodromy} of the anchor path $\delta$ is the product of the permutations $\sigma_k\dots\sigma_2\sigma_1$,  where $\rho(a_i)=\sigma_i$ is the permutation associated to the arc $a_i$ of $\alpha$.

Now we introduce notation in the statement of Theorem \ref{procedurethm}. Assume that $M_\alpha$, the 3-fold dihedral cover of $S^3$ branched along $\alpha$ corresponding to $\rho$, is a rational homology sphere. Let $\mathcal{B}_V=\{\omega_i\}_{i=1}^{2g-2}\cup \{\beta_r,\beta_l\}$ be any basis for $H_1(V-\beta;\mathbb{Z})$ consisting of embedded curves in a Seifert surface $V$ for $\alpha$, where $\beta$ is a mod 3 characteristic knot for $\alpha$, {determining $\rho$}.  Let $\delta_i$ be an anchor path for $\omega_i$, and let $\gamma_r$ and $\gamma_l$ be anchor paths for the right and left pushoffs of $\beta$ in $V$.   Let $\mu_{\delta_i}, \mu_{\gamma_r},$, and $\mu_{\gamma_l}\in D_3$ be their monodromies.  Let $c_0\in \{1,2,3\}$ be the color of the zeroth arc of $\alpha.$

 Let $\mathcal{B}=\{A_s-B_s\}_{s=1}^{2g-1}$ be the set containing the following $2g-1 $ differences of curves in $M_\alpha$:

	$$ \omega_i^j-\omega_i^k, \text{ where } \{j,k\}=\{1,2,3\}-\{\mu_{\delta_i}(c_0)\} , \text { and }$$ 
	$$\beta^j-\beta^k, \text{ where } j=\mu_{\gamma_r}(c_0) \text{ and }\{k\}=\{1,2,3\}-\{\mu_{\gamma_r}(c_0),\mu_{\gamma_l}(c_0)\}.$$
	
	 Given a simple curve $c$ on the surface $V-\beta$ we denote its positive and negative push-offs by $c^{\pm}$. 
	 We define a $(2g-1)\times (2g-1)$ matrix $M=(m_{r, s})_{r, s = 1}^{2g-1}$ as follows. Let $A_r-B_r, A_s-B_s \in \mathcal{B}$ and set  $$m_{r,s}=\text{lk}(A_r-B_r, A_s^+-B_s^-),$$ where  $\text{lk}$ denotes the linking number in $M_\alpha$.  Our main theorem states that the signature of $-M$ is the signature of the cobordism $W(\alpha, \beta)$ constructed by Cappell and Shaneson in~\cite{CS1984linking}. 
	 
 \begin{thm}\label{procedurethm}   Let $\alpha$ be a knot and $\rho:\pi_1(S^3-\alpha) \twoheadrightarrow D_3$ a surjective homomorphism given by a mod~3 characteristic knot $\beta$. Assume that $M_\alpha$, the 3-fold dihedral cover of $S^3$ branched along $\alpha$ and determined by $\rho$, is a rational homology sphere.  Let $W(\alpha,\beta)$ be the cobordism constructed in~\cite{CS1984linking} between $M_\alpha$ and the $p$-fold cyclic cover of $S^3$ branched along $\beta$. Let $M=(m_{r,s})$ be the matrix defined in the previous paragraph.  Then the signature of the 4-manifold $W(\alpha,\beta)$ is $$\sigma(W(\alpha,\beta))=-\sigma(M).$$

In particular, $\sigma(M)$ is independent of the choices of anchor paths $\delta_i$, $\gamma_r$, and $\gamma_l$ and the matrix $M$ can be used to compute the invariant $\Xi_3(\alpha)$ associated to $\rho$, using Equation~(\ref{eqXi}).
Moreover, when $M_\alpha$ is an integer homology sphere, $-M$ represents the intersection form of $W(\alpha,\beta)$.
	\end{thm}

As remarked earlier, this theorem is the non-trivial step in computing the invariant $\Xi_3(\alpha)$, since the other two terms in the formula~(\ref{eqXi}) for  $\Xi_3(\alpha)$ are determined by the Seifert forms for $\alpha$ and $\beta$, and are thus algorithmically computable from diagrams of these knots. 

In Section~\ref{examples}, we illustrate how to apply Theorem~\ref{procedurethm} to compute the signature defect associated to a singularity. We use two knots whose dihedral 3-fold covers are homeomorphic to $S^3$. Our first example is the knot $6_1$; this is the $3$-admissible knot of smallest crossing number. We use this example to illustrate a characteristic knot, anchor paths, and the associated monodromies. Our second example, the knot $8_{11}$, is a 3-admissible knot whose Seifert surface has higher genus, in order that the additional curves  $\omega_i$ and their anchor paths come into play.  In Section~\ref{octopus}, we prove Theorem \ref{procedurethm} and discuss its generalization to all odd $p$.

\section{Computing the Signature Defect}\label{examples}
 
\subsection{Overview of the procedure }

First we outline the steps for computing the signature defect $\Xi_p(\alpha)$, and, in particular, the work needed to pass from the geometric formula in \cite{kjuchukova2018dihedral} to a computation involving only diagrammatic information.  We then carry out these steps in examples.

	\begin{enumerate}
		\item Fix a diagram and Seifert surface $V$ for $\alpha$.
		\item Find a characteristic knot for $\alpha$. That is, compute the mod $p$ nullspace of the symmetrized Seifert form for $V$. Fix a primitive curve, $\beta$, in this nullspace. Choose an orientation for $\beta$.
		\item Choose a basis $\mathcal{B}_V=\{\omega_i\}_{i=1}^{2g-2}\cup \{\beta_r,\beta_l\}$ for $H_1(V-\beta;\mathbb{Z})$, where $g$ is the genus of $V$ and $\beta_r$, $\beta_l$ denote the right and left push-offs of $\beta$ in $V$.
		\item Using Theorem 1, identify the curves in the 3-fold dihedral cover of $S^3$ branched along $\alpha$ whose linking numbers contribute to the computation of $\Xi_p$.  
		\item Compute the linking numbers of these curves using the algorithm in \cite{cahnkjuchukova2016linking}. Evaluate $\sigma(W(\alpha, \beta)).$
		\end{enumerate}

\begin{ex}\label{61procedure.ex} In this example, we show $\Xi_3(6_1)=1$ using Theorem~\ref{procedurethm}. The three-coloring and the Seifert surface $V$ we use are pictured in Figure~\ref{61seifert.fig}.   Also remark that, since $6_1$ is a ribbon knot whose 3-fold dihedral cover is $S^3$, we can also conclude that $|\Xi_3(6_1)|=1$ by~\cite{cahnkjuchukova2017singbranchedcovers}. 
\begin{figure}[htbp]
	\includegraphics[width=3in]{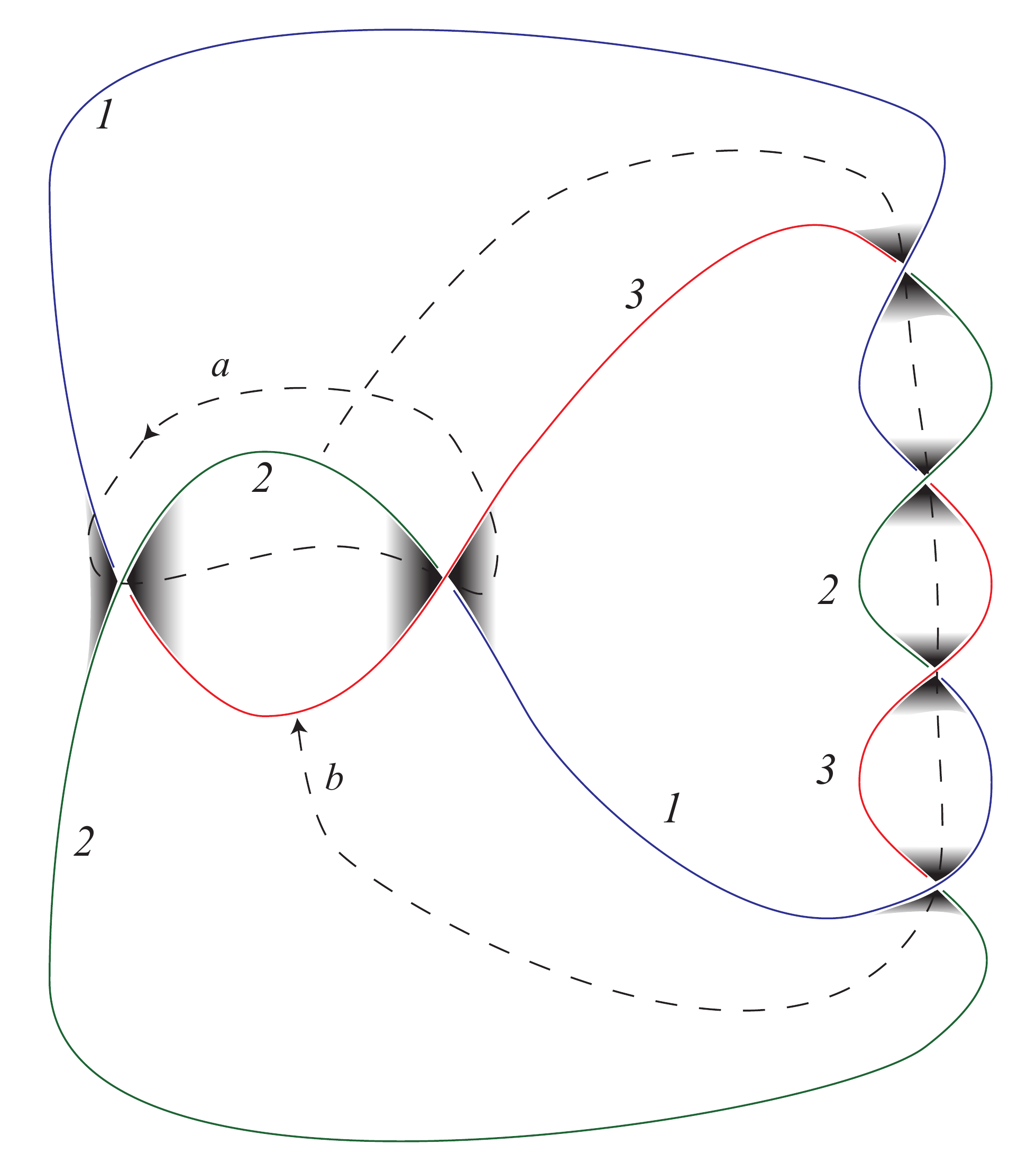}
	\caption{The knot $6_1$, a Seifert surface $V$, and a basis $\{a,b\}$ for $H_1(V;\mathbb{Z})$.}\label{61seifert.fig}
	\end{figure}
    \begin{figure}
   	\includegraphics[width=3in]{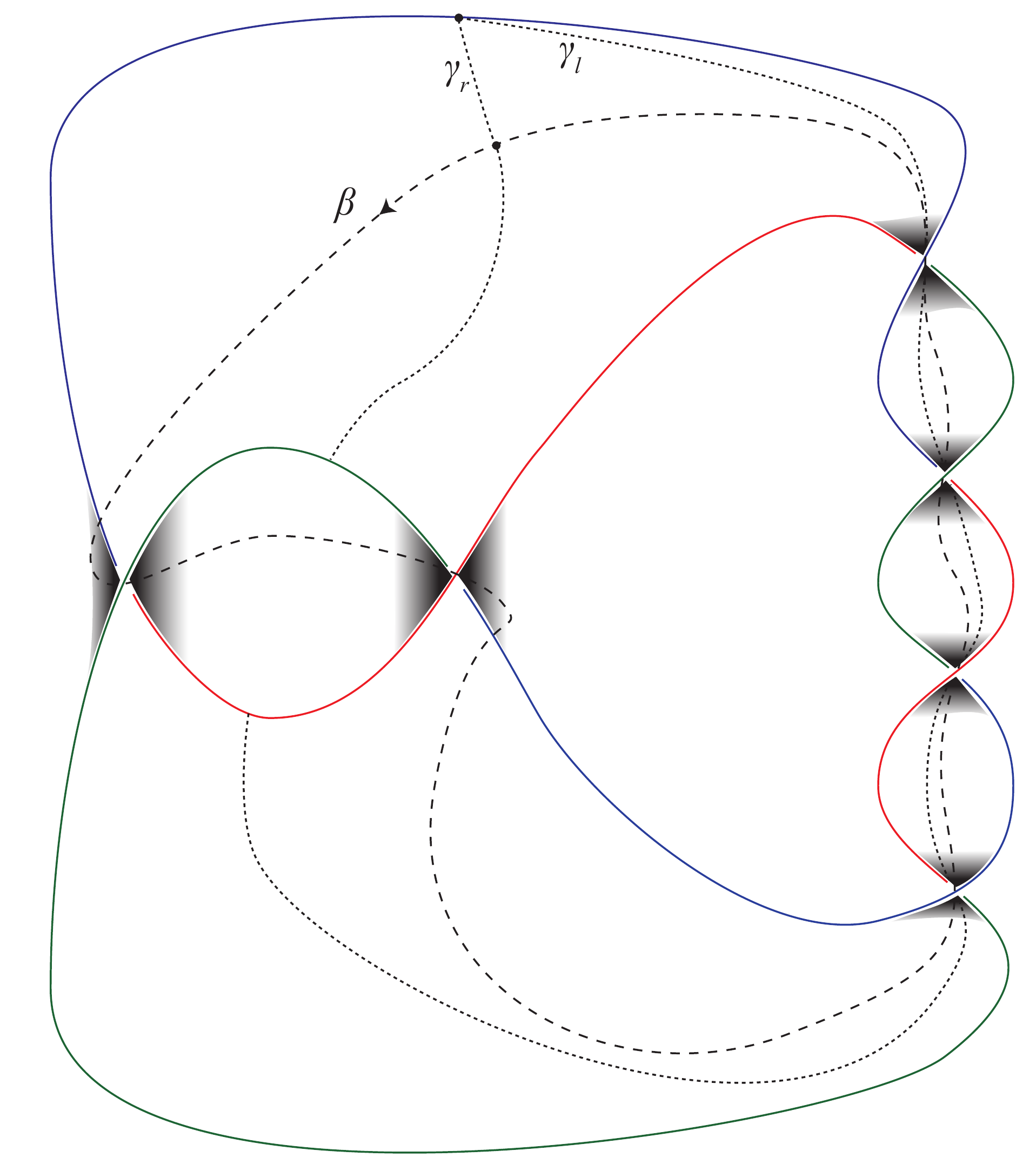}
   	\caption{The knot $6_1$ together with a characteristic knot $\beta$, and the corresponding arcs $\gamma_r$ and $\gamma_l$.}\label{61seifertwithbeta.fig}
   	\end{figure}
	
We begin by finding a mod 3 characteristic knot $\beta$ for this 3-colored diagram.  With respect to the basis $\{a,b\}$ we compute the matrix of the symmetrized linking form 
 $$[L_V]=\begin{pmatrix}-2&1\\1&4\end{pmatrix}.$$  
 
 Recall that a characteristic knot $\beta$ is one that satisfies $L_V(\beta,\omega)\equiv 0 \mod 3$ for all $\omega \in H_1(V;\mathbb{Z})$. We check that $$\begin{pmatrix}-2&1\\1&4\end{pmatrix}\begin{pmatrix}1\\-1\end{pmatrix} \equiv \begin{pmatrix}0\\0\end{pmatrix} \mod 3.$$ 
 
Hence an embedded representative of the class $a-b$ is a mod 3 characteristic knot. Moreover, since $6_1$ is a two-bridge knot, it has a single 3-coloring, up to equivalence, and therefore a single equivalence class of mod~3 characteristic knots.  Since $V$ has genus one, our basis $\mathcal{B}_V$ consists only of $\beta_r$ and $\beta_l$.  An embedded curve $\beta$, together with a choice of anchor paths $\gamma_r$ and $\gamma_l$, is shown in Figure~\ref{61seifertwithbeta.fig}.  

	\begin{figure}[htbp]
		\includegraphics[width=5in]{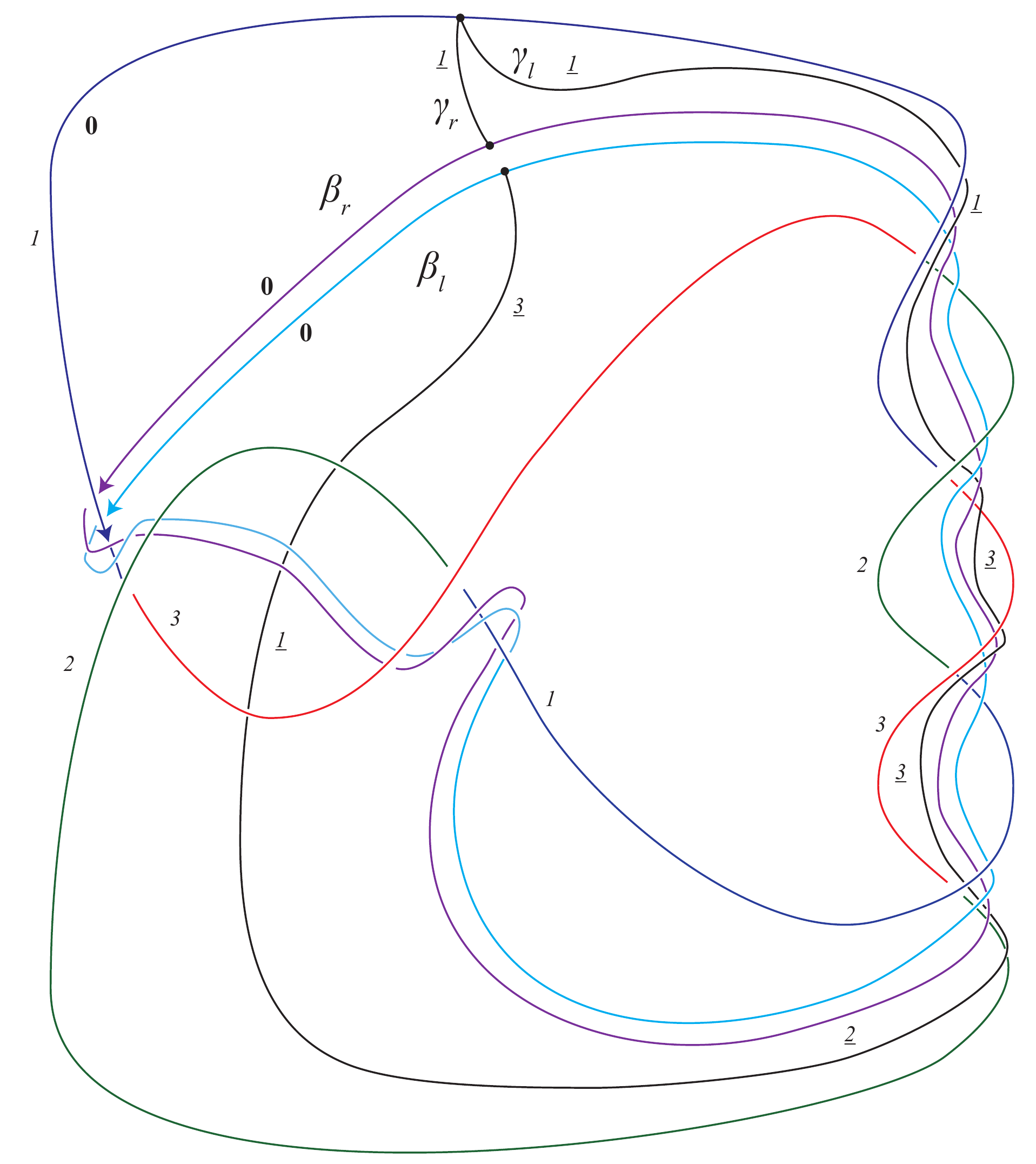}
		\caption{The link diagram for $6_1$, together with the push-offs $\beta_r$ and $\beta_l$ of the characteristic knot, and the anchor paths $\gamma_r$ and $\gamma_l$.  The arcs of $\alpha$ are marked with their colors 1, 2, or 3.  An $\underline{i}$ on an arc of $\gamma_r$ or $\gamma_l$ indicates that the lift of that arc sits in the 3-cell $e^3_i$.}
		\label{6-1numbereddiagram.fig}
		\end{figure}

We use the algorithm in \cite{cahnkjuchukova2016linking} to compute linking numbers of the lifts $\beta_r^j$ and $\beta_l^k$, $j,k=1,2,3$.  Details of this computation are given in the Appendix. The $jk$-entry of the matrix $$\begin{pmatrix}-1&0&1\\0&0&0\\1&0&-1\end{pmatrix}$$  

is the linking number of $\beta_r^j$ and $\beta_l^k$.  As $\beta_r$ and $\beta_l$ are parallel curves in $V$, we may also view this as the matrix of linking numbers of $\beta^j$ with $\beta^k$, in which case the diagonal entries may be interpreted as self-linking numbers.

Next we compute the monodromies of $\gamma_r$ and $\gamma_l$, in order to determine which of the above linking numbers appear in the formula for $\Xi_3(6_1)$:
$$\mu_{\gamma_r}=Id$$
$$\mu_{\gamma_l}=(23)(13)(12)(23)(12)(13)=(123)$$

The zeroth arc of $\alpha$ is colored $c_0=1$.  Hence $\mu_{\gamma_r}(c_0)=1$ and $\mu_{\gamma_l}(c_0)=3$.  One can see this in the link diagram in Figure~\ref{6-1numbereddiagram.fig}; the underlined number $\underline{i}$ on an arc of $\gamma_r$ or $\gamma_l$ indicates that the relevant lift of that arc sits in the 3-cell $e^3_i$.  These cells change from one arc of $\gamma_r$ or $\gamma_l$ to the next according to $\mu_{\gamma_r}$ and $\mu_{\gamma_l}$.

By Theorem~\ref{procedurethm}, the signature $\sigma(W(\alpha,\beta))$ is the signature of the $1\times 1$ matrix whose entry is the linking of $\beta^1-\beta^2$ with itself.  Using the linking numbers given in the matrix above, we see that $\text{lk}(\beta^1-\beta^2,\beta^{1,+}-\beta^{2,-})=-1$.  Hence $\sigma(W(\alpha,\beta))=-\sigma((-1))=1$.  Since $\beta$ is an unknot with zero self-linking, the other terms in the formula for $\Xi_3(6_1)$ vanish and we conclude that $\Xi_3(6_1)=1$.

\end{ex}

\begin{ex}\label{811procedure.ex} In this example, we show $\Xi_3(8_{11})=3$. This answer is independently confirmed in~\cite{cahnkju2018genus} using a trisection diagram of an irregular dihedral branched cover of $S^4$. We remark that the technique used in~\cite{cahnkju2018genus} to evaluate the invariant $\Xi_3(\alpha)$ is rather less computationally onerous. However, this alternative method can only be applied when given a Fox colored triplane diagram of a branching set in $S^4$ with a singularity of type $\alpha$. By contrast, the procedure used here is inherently 3-dimensional; it only uses a colored diagram of $\alpha$.

	Let $V$ be the Seifert surface for $8_{11}$ given by checkerboard coloring the diagram in Figure~\ref{8_11pic.fig}.
	 With respect to the basis $\{A,B,\gamma,\beta\}$ of $H_1(V;\mathbb{Z})$ in Figure~\ref{8_11pic.fig}, where $\gamma=\gamma_r\cdot\gamma_l$, we find that the matrix of the symmetrized Seifert form is
	$$[L_V]=\begin{pmatrix} 
	2&-1&0&0\\ 
	-1&2&-1&0\\
	0&-1&-2&-3\\
	0&0&-3&0
	\end{pmatrix}$$
	
  We verify that
	$$\begin{pmatrix} 
	2&-1&0&0\\ 
	-1&2&-1&0\\
	0&-1&-2&-3\\
	0&0&-3&0
	\end{pmatrix}\begin{pmatrix}0\\0\\0\\1\end{pmatrix}\equiv \begin{pmatrix}0\\0\\0\\0\end{pmatrix} \mod 3,$$ so the curve $\beta$ is a characteristic knot. As in Example~1, any mod~3 characteristic knot for $8_{11}$ determines its unique 3-coloring.
	
	A basis $\mathcal{B}_V$ for $H_1(V-\beta;\mathbb{Z})$ is $\{A,B,\beta_r,\beta_l\}$.  We again use the algorithm in \cite{cahnkjuchukova2016linking} to find all linking numbers of lifts of curves in $\mathcal{B}_V$.  These linking numbers are displayed in Table ~\ref{8_11linkingnums.tab}.   We use superscript $\pm$ to denote positive and negative push-offs of curves on $V$, as well as their lifts.

	\begin{table}[h]
	\begin{tabular}{c|c|c|c|}
		& $A^+$& $B^+$ &$\beta^+$\\ \hline
	$A$&$\begin{pmatrix} 1&0&0\\0&0&1\\0&1&0\end{pmatrix}$&$\begin{pmatrix} 0&0&0\\0&0&0\\0&0&0\end{pmatrix}$& 	$\begin{pmatrix} 0&0&0\\-1&0&1\\1&0&-1\end{pmatrix}$\\ \hline
	$B$&$\begin{pmatrix} -1&0&0\\0&-1&0\\0&0&-1\end{pmatrix}$&$\begin{pmatrix} 1&0&0\\0&0&1\\0&1&0\end{pmatrix}$&$\begin{pmatrix} 0&0&0\\-1&0&1\\1&0&-1\end{pmatrix}$\\\hline
	$\beta$& $\begin{pmatrix} 0&-1&1\\0&0&0\\0&1&-1\end{pmatrix}$&$\begin{pmatrix} 0&-1&1\\0&0&0\\0&1&-1\end{pmatrix}$&$\begin{pmatrix} -3&0&3\\0&0&0\\3&0&-3\end{pmatrix}$\\\hline
	\end{tabular}	
	\caption{Linking numbers of lifts of the curves $\{A,B,\beta\}$ and their push-offs on $V$.}\label{8_11linkingnums.tab}
	\end{table}
	
	\begin{table}[h]
	\begin{tabular}{c|c|c|c|}
	&$A^{2,+}-A^{3,-}$&$B^{2,+}-B^{3,-}$&$\beta^{1,+}-\beta^{2,-}$\\\hline
	$A^2-A^3$&-2&-1&-2\\\hline
	$B^2-B^3$&-1&-2&-2\\\hline
	
	$\beta^1-\beta^2$&-2&-2&-3\\\hline
	\end{tabular}
	\caption{ The matrix $M$ of linking numbers of curves in $\mathcal{B}$.  The intersection form of $W(\alpha,\beta)$ is given by the matrix $-M$.}\label{8_11matrix.tab}
	\end{table}

The matrix $M=\begin{pmatrix}-2&-1&-2\\ -1&-2&-2\\-2&-2&-3\end{pmatrix}$ of Theorem~\ref{procedurethm}, shown in Table ~\ref{8_11matrix.tab}, has signature $-3$.  Since $\beta$ is an unknot and $L_V(\beta,\beta)=0$, we have $\Xi_3(8_{11})= \sigma (W(\alpha, \beta)) =-\sigma(M)=3$.
\begin{figure}[htbp]
	\includegraphics[width=6in]{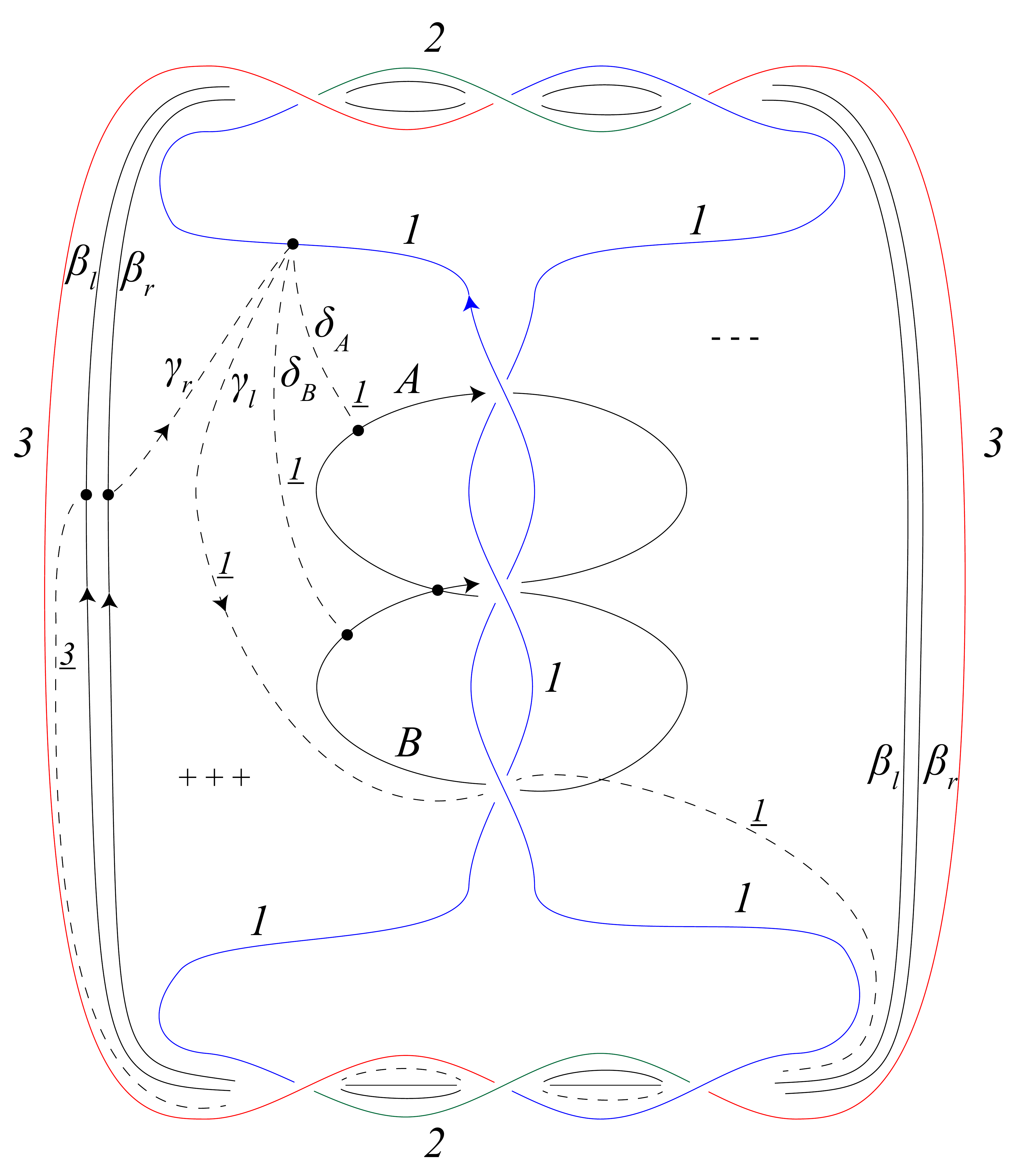}
	\caption{The knot $8_{11}$, shown in a 9-crossing diagram, with a choice of characteristic knot and anchor paths.}\label{8_11pic.fig}
\end{figure} 
\end{ex}

		\section{Proof of Theorem~\ref{procedurethm} and the Cappell-Shaneson construction}\label{octopus}
		
Before proving Theorem~\ref{procedurethm} we briefly review the Cappell-Shaneson construction of $M_\alpha$, the irregular dihedral cover of $S^3$ branched along $\alpha$, and a cobordism $W(\alpha,\beta)$ between $M_\alpha$ and $M_\beta$, the $p$-fold cyclic cover of $S^3$ along $\beta$, where $\beta$ is a mod $p$ characteristic knot for $\alpha$.  We again focus on the case $p=3$, but our combinatorial procedure can be generalized to all odd $p$, as discussed at the end of this section.
	
\subsection{The Cappell-Shaneson Construction of the Irregular Dihedral Cover}\label{dih-con}   
In~\cite{CS1984linking}, Cappell and Shaneson construct the irregular $p$-fold dihedral cover of $S^3$ branched along $\alpha$ from the $p$-fold cyclic cover of $S^3$ branched along a characteristic knot $\beta$. In this paragraph, we give an informal overview of their construction.  Precise details will be provided later, as needed.  Roughly speaking, one begins with the $p$-fold cyclic cover $M_\beta$, and considers the lifts to this cover of a Seifert surface for $\alpha$, $V\supset \beta$. Remove from $M_\beta$  a neighborhood $J$ of the union of the preimages of $V$ to obtain a 3-manifold with boundary $\partial J$. Now identify points on that boundary via an involution $\bar{h}$, defined below.  The resulting closed manifold $M_\alpha$ is the $p$-fold irregular dihedral cover of $S^3$ branched along $\alpha$.  The surface $S:=\bar{h}(\partial J)$ sits inside this covering space, and has boundary equal to the index-one lift of $\alpha$.  The index-two lift of $\alpha$ is an embedded curve on $S$.  In order to compute the signature $\sigma(W(\alpha,\beta))$ in Theorem~\ref{procedurethm}, we must compute a matrix of linking numbers of certain elements of $H_1(S;\mathbb{Z})$; namely, a basis for the kernel of the map $i_*:H_1(S; \mathbb{Z})\rightarrow H_1(\bar{h}(J);\mathbb{Z})$.

\begin{figure}[htbp]
	\begin{subfigure}[b]{.4\textwidth}
	\includegraphics[width=2.8in]{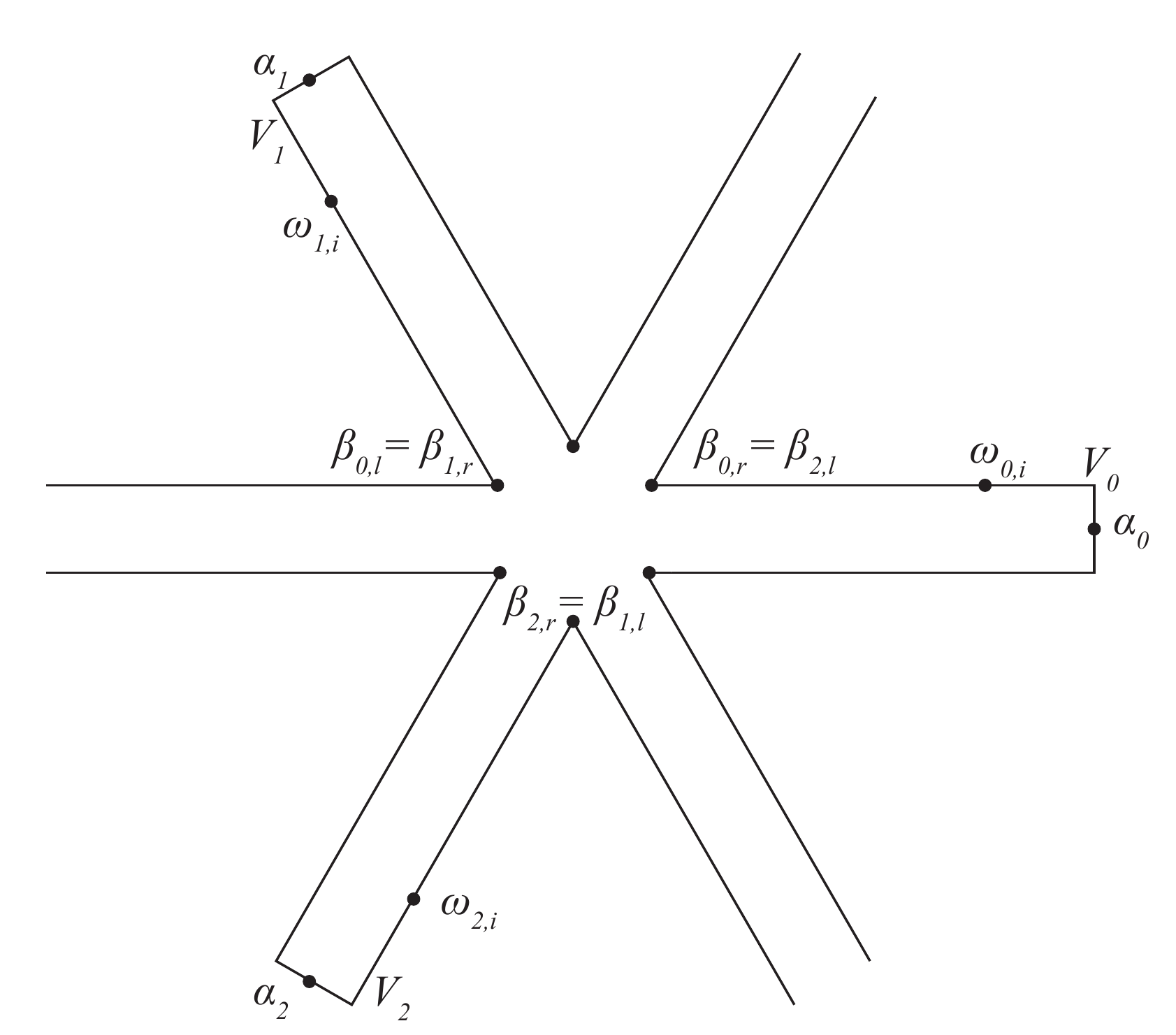}
	\end{subfigure} \quad\quad\quad\raisebox{1in}{$\xrightarrow{\bar{h}}$}\quad
	\begin{subfigure}[b]{.4\textwidth}
	\raisebox{.75in}{\includegraphics[width=2.8in]{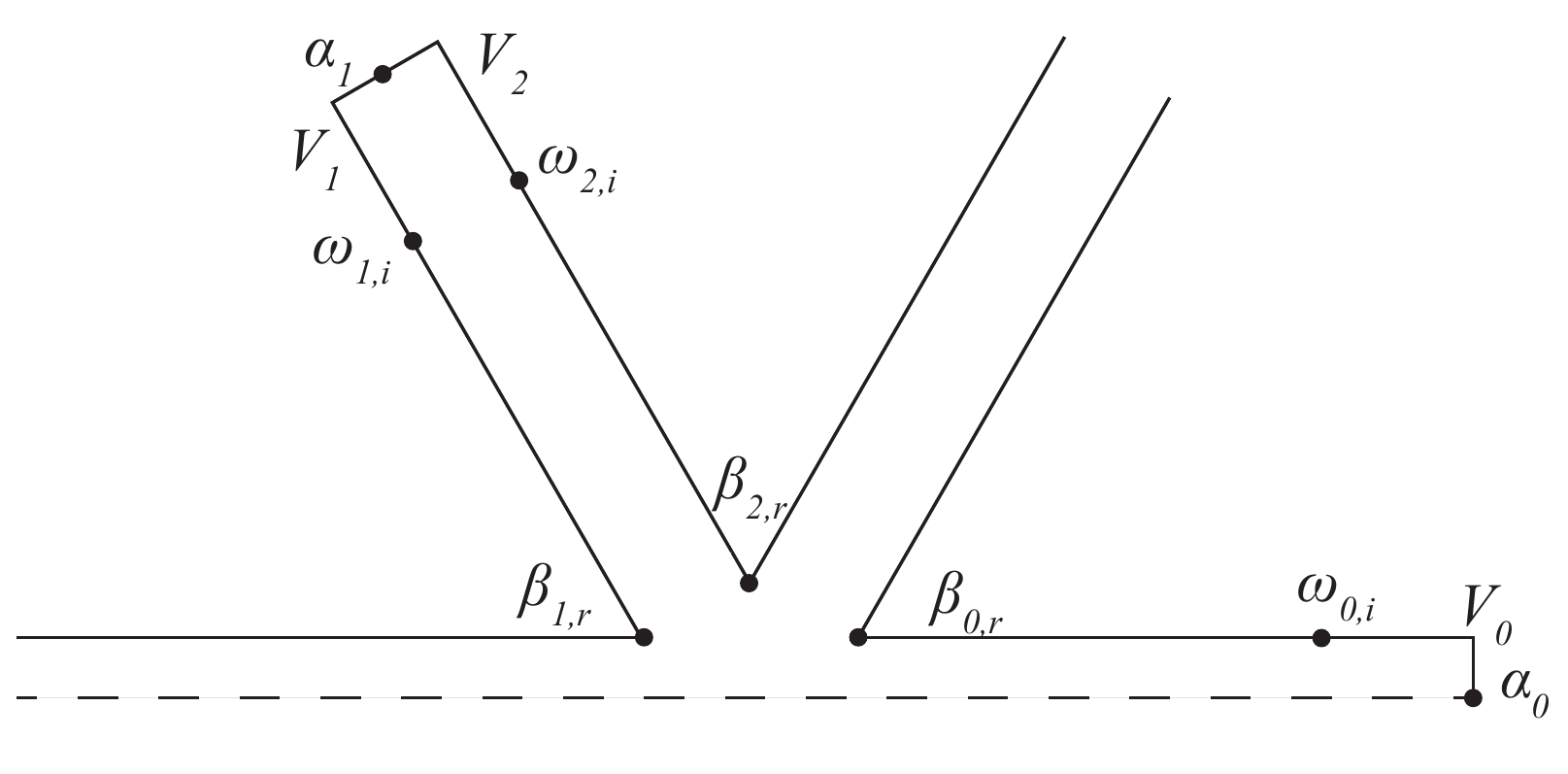}}
	\end{subfigure}
	\caption{The image of $J$ under the involution $\bar{h}$.} \label{halfstar.fig}
	\end{figure}

Now we set $p=3$, describe the construction in detail, and introduce the necessary notation.   Let $f:M_\beta\rightarrow S^3$ be a $3$-fold cyclic covering map branched along $\beta$.  By the construction sketched out above, we know that $M_{\alpha}$ can be obtained from $M_\beta$ as follows. Let $J=f^{-1}(V\times[-1,1])\subset M_\beta$. Let $h:V\times[-1,1]\rightarrow V\times[-1,1]$ be given by $h(x,t)=h(x,-t)$.  Let $\bar{h}:\partial J\rightarrow \partial J$ be the lift of $h$ to $M_\beta$ restricted to $\partial J$; in the schematic in Figure ~\ref{halfstar.fig}, $\bar{h}$ is a reflection about the horizontal line.  Cappell and Shaneson show that $M_{\alpha}$ is homeomorphic to $(M_\beta-\mathring{J})/\{\bar{h}(x){\sim}x, x\in \partial J\}$, and that the mapping cone $W(\alpha,\beta)$ of $\bar{h}$ is a cobordism from the $3$-fold cyclic cover $M_\beta$ to the irregular 3-fold dihedral cover $M_{\alpha}$.  The surface $S=\bar{h}(\partial J)$ is embedded in $M_{\alpha}$, and has one boundary component $\alpha_0$, the index-one lift of $\alpha$. 

\begin{figure}[htbp]\includegraphics[width=3in]{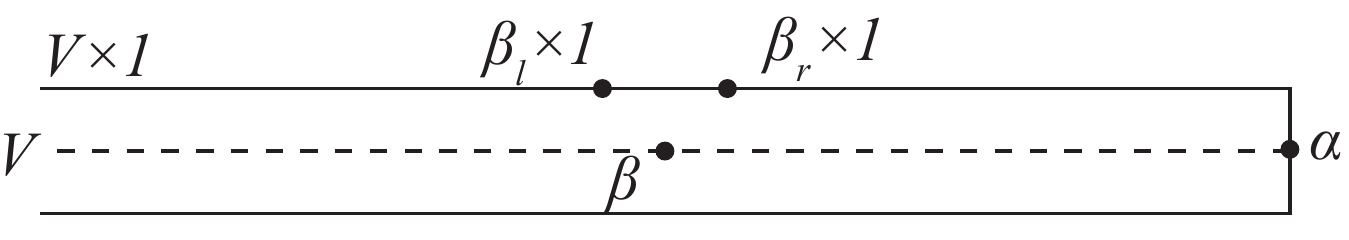}
	\caption{$V$ is a Seifert surface for $\alpha$.  The needed push-offs of $\beta$ in $V\times[-1,1]$ are shown above.}\label{VtimesI.fig}
	\end{figure}

Let $V-\beta$ denote the surface $V$ cut along $\beta$, which we obtain by removing a thin annulus between the right and left push-offs $\beta_r$ and $\beta_l$ of $\beta$ in $V$ (note that $\beta$ is oriented). The surface $S$ above can be obtained by gluing together three copies of $V-\beta$ as follows.  There are three lifts of $(V-\beta)\times 1$ in $M_\beta$, which we label $V_0$, $V_1$, and $V_2$, according to the action of the deck transformation group.  Let $\alpha_0$, $\alpha_1$, and $\alpha_2$ denote the corresponding lifts of $\alpha$. Each $V_i$ contains lifts of the curves $\beta_r\times 1$ and $\beta_l\times 1$, and we denote these by $\beta_{i,r}$ and $\beta_{i,l}$.  See Figures \ref{VtimesI.fig} and \ref{halfstar.fig}.

 From Figure \ref{halfstar.fig}, we can read off the boundaries of the surfaces $V_i$:

$$\partial V_0= \alpha_0 +\beta_{0,r}-\beta_{0,l}$$
$$\partial V_1=\alpha_1 + \beta_{1,r}-\beta_{1,l} $$
$$\partial V_2=\alpha_2 + \beta_{2,r}-\beta_{2,l}.$$

\begin{figure}[htbp]
	\includegraphics[width=3in]{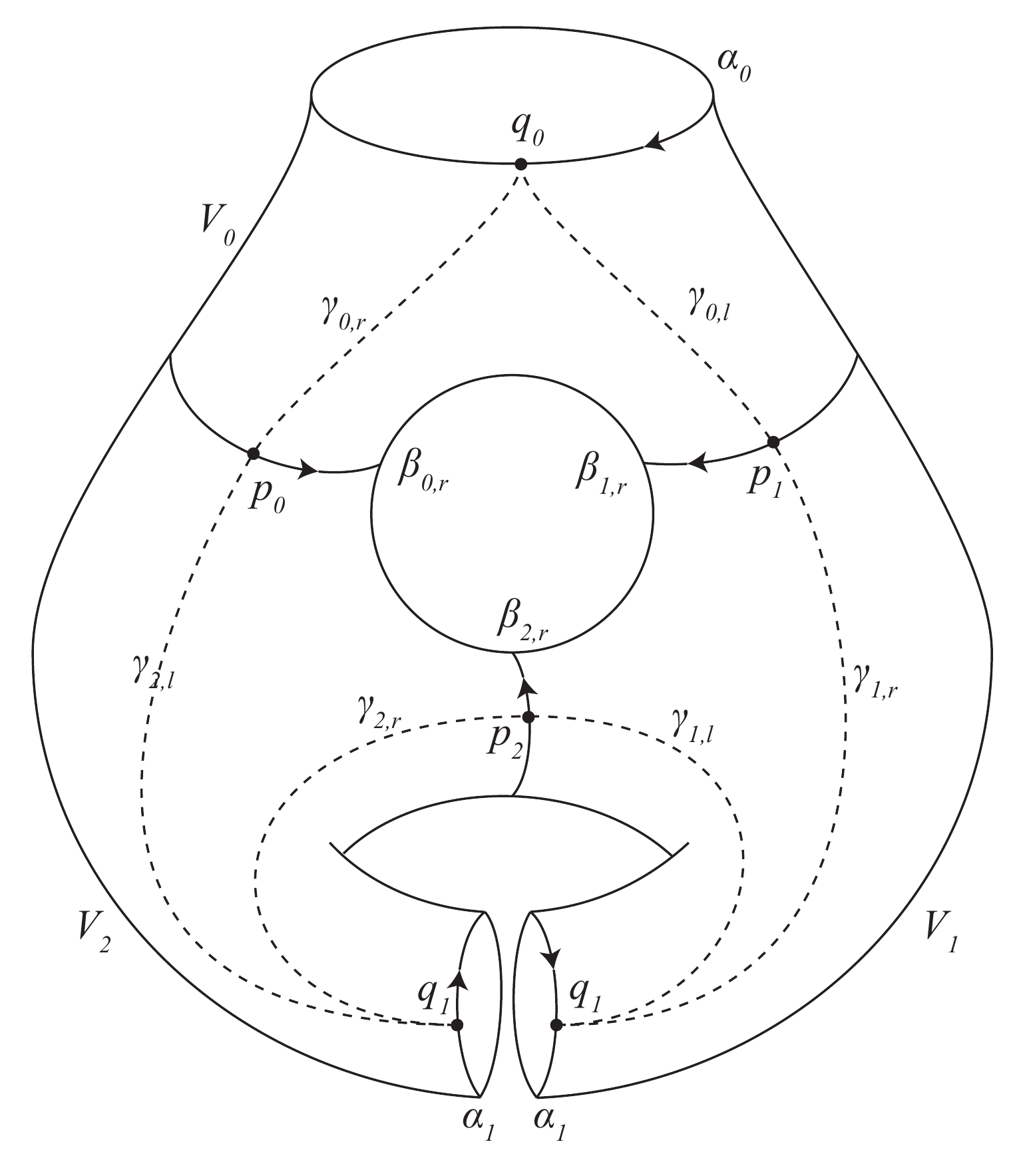}
	\caption{The surface $S$, together with the index-one and index-two branch cuves $\alpha_0$ and $\alpha_1$, the three lifts of $\beta$, and the three lifts of $\gamma_r$ and $\gamma_l$.  The two copies of $\alpha_1$ are identified according to the orientations pictured above.}
\label{halfoctopus.fig}
	\end{figure}

Now we construct $S$ by gluing together $V_0$, $V_1$, and $V_2$ using the following identifications: $\beta_{0,l}$ is identified with $\beta_{1,r}$, $\beta_{1,l}$ is identified with $\beta_{2,r}$, and $\beta_{2,l}$ is identified with $\beta_{0,r}$.  In addition $\alpha_1$ and $\alpha_2$ are identified.  The index-one and index-two branch curves are $\alpha_0=\partial S$ and $\alpha_1$ respectively.  Note that $\beta_{0,r}$ and $\beta_{1,r}$ are homologous in $S$, as they cobound $V_0$ together with $\alpha_0$.  The surface $S$, constructed using these identifications, is pictured in Figure~\ref{halfoctopus.fig}, in the case where $V$ has genus one and each $V_i$ is a pair of pants.  This is in fact the case in our first example, where $\alpha$ is the knot $6_1$.  In general the genus of $V_i$ is one less than the genus of $V$.

\subsection{Proof of Theorem \ref{procedurethm}}  Let $W=W(\alpha,\beta)$ denote the Cappell-Shaneson cobordism, described above, between $M_\alpha$, the $p$-fold dihedral branched cover of $\alpha$, and $M_\beta$, the $p$-fold cyclic branched cover of the characteristic knot $\beta$. We will show that   $\sigma(W) = - \sigma(M)$, where  $M$ is the matrix of linking numbers between the set $\mathcal{B}$ of curves in $M_\alpha$ from the theorem statement. The argument is standard: identify relative classes in $H_2(W, M_\alpha; \mathbb{Z})$ with curves in $\mathcal{B}$, and then show that the intersection numbers of two 2-dimensional classes equal the linking numbers of the corresponding curves.  One remark we make right away is that, to calculate the intersection number, we actually use the linking numbers in $\overline{M}_\alpha$, that is, in $M_\alpha$ with the orientation reversed.  This is a matter of convention. We adopted the convention used in~\cite{cahnkju2018genus} which means that $M_\alpha$ considered as a cover of $S^3$ branched along $\alpha$ has the opposite orientation than the one it inherits as a boundary component of the 4-manifold $W$. Hence, a minus sign appears in the signature formula, $\sigma(W)=-\sigma(M)$, since $M$ is a matrix of linking numbers in in $M_\alpha$. Using the opposite convention would amount to replacing $\alpha$ by its mirror, in which case the orientation reversal would be unnecessary.

We have that $M_\alpha$ is a rational homology sphere by assumption  and $M_\beta$ is a rational homology sphere because it is the 3-fold branched cover of $S^3$ along a knot. Thus, we can identify $H_2(W; \mathbb{Z})$ with its image in any of the relative groups  $H_2(W, M_\alpha; \mathbb{Z})$,  $H_2(W, M_\beta; \mathbb{Z})$ or  $H_2(W, M_\alpha\cup M_\beta; \mathbb{Z})$, since the inclusion map in each of the relevant long exact sequences is injective. We will work with the image $\iota_\ast(H_2(W; \mathbb{Z}))\subset H_2(W, M_\alpha; \mathbb{Z})$, denoted by $H$ from now on, since we happen to have a basis for $H_2(W, M_\alpha;\mathbb{Z})$ on hand. Indeed, by ~\cite[Equation 2.23]{kjuchukova2018dihedral}, $H_2(W, M_\alpha;\mathbb{Z})$ is {free and} isomorphic to $\ker i_*$, where $i_*:H_1(S; \mathbb{Z})\rightarrow H_1(\bar{h}(J);\mathbb{Z})$.

Corollary~2.4 of~\cite{kjuchukova2018dihedral} describes a basis for $\ker i_*\subset H_1(M_\alpha; \mathbb{Z})$, where each element in this kernel is identified with a relative cycle in $H_2(W, M_\alpha;\mathbb{Z})$ via the exact sequence of the pair.   A relative cycle in $H_2(W,M_\alpha;\mathbb{Z})$ lies in  $H=\iota_\ast(H_2(W; \mathbb{Z}))$ if and only if its boundary can be capped off by {an oriented} surface in $M_\alpha$. Given two classes in  $H$, we will compute their intersection number in terms of linking numbers between their boundaries, which we will describe in terms of $\ker i_*$.

 The remaining ingredient is to characterize the curves in a basis for $\ker i_*$, which lie in the dihedral cover $M_\alpha$,  using only diagrammatic data about the branching set $\alpha$ in $S^3$. The curves in $\ker i_*$ project under the branched covering map $M_\alpha\to S^3$ to curves in $S^3-\alpha$ and can be described in terms of a basis for $H_1(V-\beta; \mathbb{Z})$~\cite{kjuchukova2018dihedral}. Each curve in $H_1(V-\beta; \mathbb{Z})$ is covered by 3 disjoint circles in $M_\alpha$ as seen in the Cappell-Shaneson construction~\cite{CS1984linking}. We use anchor paths and their monodromies to give a combinatorial description of $\ker i_*$ relying solely on diagrammatic information in $S^3$. We conclude that the curves in $\mathcal{B}$ represent a basis for  $\text{ker }i_*$ and, equivalently, for $H_2(W,M_\alpha;\mathbb{Z})$.

		Let $q$ be a point on $\alpha$, and let $p$ be a point on $\beta$. Let $\{\omega_1,\dots,\omega_{2g-2}\}\cup \{\beta_r,\beta_l\}$ be a basis for $H_1(V-\beta;\mathbb{Z})$, where $g$ is the genus of $V$.  Each curve $\omega_i$ in $V-\beta$ has three lifts $\omega_{0,i}$, $\omega_{1,i}$, and $\omega_{2,i}$ to $S\subset M_{\alpha}$. From Figure \ref{halfstar.fig}, it is evident that the differences of curves $\omega_{1,i}-\omega_{2,i}$, together with one of either  $\beta_{0,r}- \beta_{2,r}$ or  $\beta_{1,r}- \beta_{2,r}$, form a basis $\mathcal{B}$ for $\ker i_*$. 
		We refer the reader to Section~\ref{algoverview.sec} for the notation and definitions used in Theorem ~\ref{procedurethm}, namely the lifts $\omega^j_i$, $\beta^j_r$, $\beta^j_l$, $\beta^j$, and the definitions of the anchor paths $\gamma_r$, $\gamma_l$, $\delta_i$, and their monodromies.

Now we use anchor paths in a diagram of $\alpha$ to identify which of the three lifts $\beta^1_r$, $\beta^2_r$, and $\beta^3_r$ are $\beta_{0,r}$, $\beta_{1,r}$, and $\beta_{2,r}$.  The lift $\gamma_{0,r}$ of the anchor path $\gamma_r$ beginning at $q_0$ has its initial endpoint in the 3-cell $e^3_{c_0}$, on the index-one lift $\alpha_0$ of $\alpha$.  Looking at $\gamma_r$ in a diagram of $\alpha$ (see, e.g., Figure ~\ref{6-1numbereddiagram.fig}), we see $\gamma_{0,r}$ has its final endpoint on the lift $\beta^j_r$ of $\beta_r$ that lies in the 3-cell $e^3_j$, where $j=\mu_{\gamma_r(c_0)}$.   On the other hand, from Figure~\ref{halfoctopus.fig}, we see that the final endpoint of $\gamma_{0,r}$ lies on the lift $\beta_{0,r}$ of $\beta.$ Hence $\beta_{0,r}=\beta^j_r$ with $j=\mu_{\gamma_r(c_0)}$. Similarly, the lift $\gamma_{0,l}$ of $\gamma_l$ beginning at $q_0$ has its initial point in the 3-cell $e^3_{c_0}$, on the index-one lift $\alpha_0$ of $\alpha$, and its endpoint on the lift $\beta^k_l$ of $\beta_l$ that lies in the 3-cell $e^3_k$, where $k=\mu_{\gamma_l(c_0)}$.  We also see from Figure~\ref{halfoctopus.fig} that the endpoint of $\gamma_{0,l}$ lies on the lift $\beta_{1,r}$ (which is identified with $\beta_{0,l}$) of $\beta_r$, so $\beta_{1,r}=\beta^k_l$ with $k=\mu_{\gamma_l(c_0)}$.  Hence the basis element $\beta_{0,r}-\beta_{2,r}$ of $\ker i_*$ is 
$$\beta^j_r-\beta^k_l, \text{ where } j=\mu_{\gamma_r}(c_0) \text{ and }\{k\}=\{1,2,3\}-\{\mu_{\gamma_r}(c_0),\mu_{\gamma_l}(c_0)\},$$
as stated in Theorem ~\ref{procedurethm}.  Note that we may instead use $\beta^j-\beta^k$, as $\beta_r$ and $\beta_l$ are simply push-offs of $\beta$ in $V$.

\begin{figure}\includegraphics[width=3in]{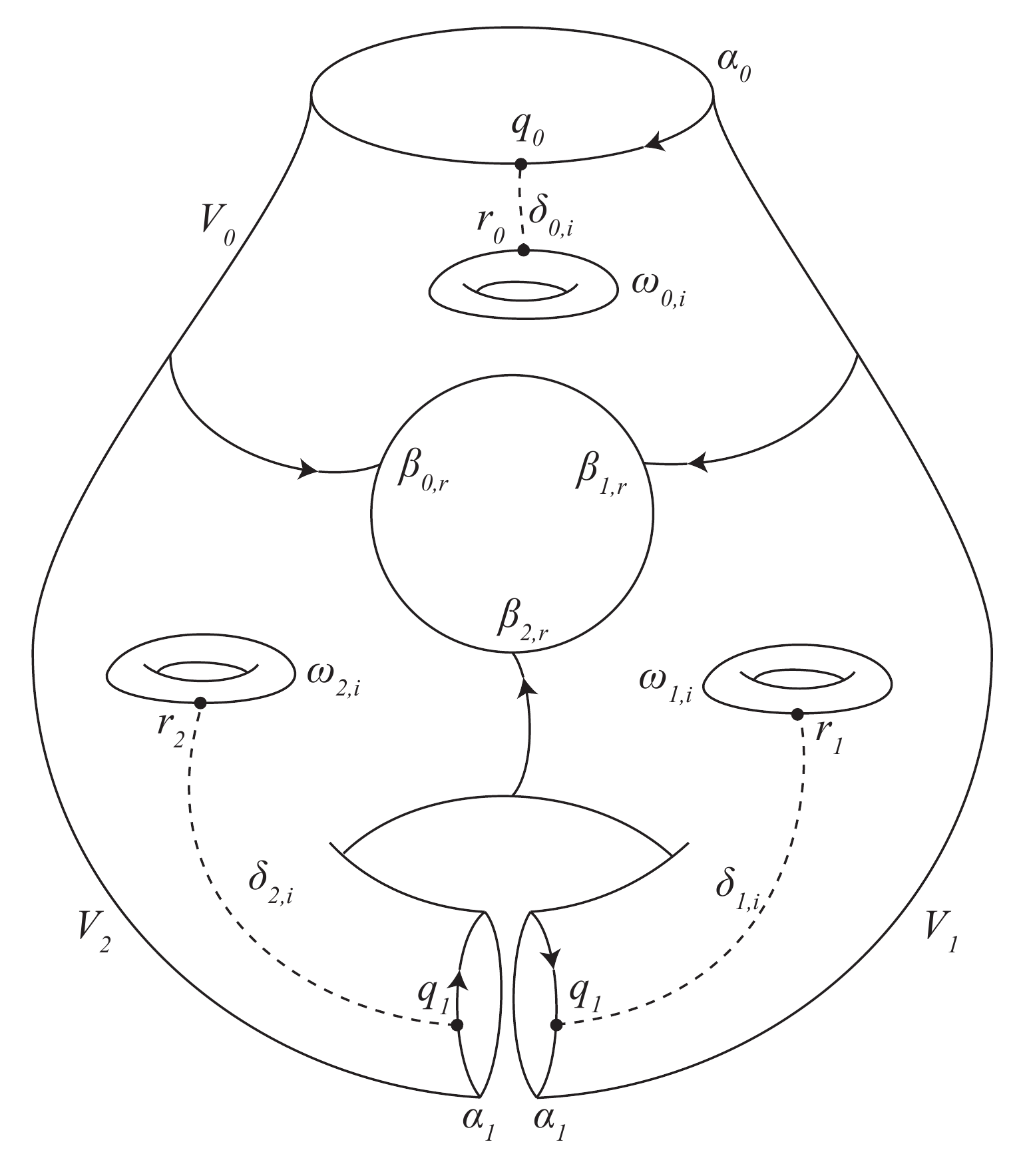}
	 \caption{Lifts to the surface $S$ of an anchor path for $\omega_i$. The two copies of $\alpha_1$ are identified according to the orientations pictured above.}
	 \label{halfoctopuswithomega.fig}
	 \end{figure}

Next we need to identify the basis element $\omega_{1,i}-\omega_{2,i}$ using diagrammatic information.  The lift $\delta_{i,0}$ of $\delta_i$ beginning at $q_0$ has its initial endpoint in the 3-cell $e^3_{c_0}$, on the index-one lift $\alpha_0$ of $\alpha$.  Looking at $\delta_i$ in a diagram of $\alpha$, we see that $\delta_{i,0}$ has its final endpoint on the lift $\omega^j_i$ where $j=\mu_{\delta_i}(c_0)$.  On the other hand, from Figure~\ref{halfoctopuswithomega.fig}, we see that the final endpoint of $\delta_{0,i}$ lies on the lift $\omega_{0,i}$ of $\omega_i$.  Hence $\omega_{0,i}=\omega^j_i$ where $j=\mu_{\delta_i}(c_0)$.  Therefore the basis element $\omega_{1,i}-\omega_{2,i}$ is, up to sign, 

$$\omega_i^j-\omega_i^k, \text{ where } \{j,k\}=\{1,2,3\}-\{\mu_{\delta_i}(c_0)\},$$
as stated in Theorem 1.

So far, we have seen that the curves in $\mathcal{B}$ can be identified with a basis for $H_2(W, M_\alpha; \mathbb{Z})\supset H$, and we wish to calculate the intersection numbers of classes in $H$. Given any $A_r-B_r\in \mathcal{B}$, we have $ A_r-B_r = \partial C_r$ where $[C_r]\in H_2(W,M_\alpha;\mathbb{Z})$ can be represented by a cylinder~\cite[Corollary 2.4]{kjuchukova2018dihedral}. 

First suppose $M_\alpha$ is an integer homology sphere. In this case there exist Seifert surfaces $\Sigma_r=\partial A_r$ and $\Sigma'_r=\partial B_r$.  The classes represented by the closed surfaces 

$$\hat{C}_r:= C_r \cup \Sigma_r \cup \Sigma_r'$$

form a basis for $H$.  Note that $A_r-B_r$ and $A_s-B_s$ do not, in general, bound disjoint cycles in $H_2(W,M_\alpha; \mathbb{Z})$,  but $A_r-B_r$ and $A_s^+-B_s^-$ bound disjoint cylinders for any choice of $r, s$. This is because $A_r-B_r$ bounds a cylinder contained entirely in $\bar{h}(J)$, in the notation of Section~\ref{dih-con}. To give an example of such a cylinder, set $A_r-B_r=\omega_{1, i}-\omega_{2, i}$, pictured as a pair of points in Figure~\ref{halfstar.fig}. In this figure, the cylinder $S^1\times I$ would be depicted as a line segment, properly embedded in  $\bar{h}(J)$, connecting $\omega_{1, i}$ to $\omega_{2, i}$. The same holds for all differences of curves in $\mathcal{B}$. On the other hand, the push-off $A_s^+-B_s^-$ bounds a cylinder which can be made disjoint from all of $\bar{h}(J)$. Therefore, the intersection number of $[\hat{C}_r]$ and $[\hat{C}_s]$ is equal to the linking number of $A_r-B_r$ with $A_s^+-B_s^-$ in $\overline{M}_\alpha$, where again, we reverse the orientation on $M_\alpha$ for the reason previously explained.

 Let $M=(m_{r,s})$ be the matrix with entries $$m_{r,s}=\text{lk}(A_r-B_r, A_s^+-B_s^-)$$ where  $A_r-B_r, A_s-B_s$  run through the elements of $\mathcal{B}$. If $M_\alpha$ is an integer homology sphere, $\mathcal{B}$ is in fact a basis for $H$, and the matrix $-M$ represents the intersection form of $W$. In particular, $\sigma(W)=-\sigma(M)$ as claimed. If $M_\alpha$ is only a rational homology sphere, we can express a basis $\mathcal{B}'$ for $H$ in terms of the $C_r$. Given two elements $\Sigma_i a_i C_{r_i}, \Sigma_j b_j C_{s_j} \in H$, their intersection number is $$- \Sigma_{i,j} a_ib_j \text{lk}(A_{r_i}-B_{r_i}, A^+_{s_j}-B^-_{s_j})=-\Sigma_{i,j}a_ib_jm_{{r_i},{s_j}}.$$ Since $H$ and $H_2(W, M_\alpha)$ have the same rank, the intersection matrix with respect to $\mathcal{B}'$ and  $-M$ are congruent over $\mathbb{Q}$, so they have the same signature.~\qedsymbol

	 \subsection{Computation for other values of $p$}
	 We briefly explain how these methods generalize to an algorithm for computing $\Xi_p(\alpha)$ for all odd $p$. The schematic in Figure~\ref{halfstar.fig} allows one to read off a basis for $\ker i_*$ without the hypothesis $p=3$. (Of course, in the general case, one uses the $p$-fold cyclic cover of $V\times [-1,1]$ branched along $\beta$, rather than the 3-fold cover.)  Anchor paths can be used in the same fashion to distinguish between the $p$ lifts of $\beta$ and the $p$ lifts of each $\omega_i$ in the diagram for $\alpha=\partial V$.  The remaining task is the computation of the linking numbers of these lifts in the $p$-fold dihedral cover of $S^3$ branched along $\alpha$. This algorithm carried out in \cite{cahnkjuchukova2016linking} in the case $p=3$ generalizes to any odd $p$.

\section*{Appendix: Calculating linking numbers in branched covers}
\subsection{Making lists.}
We explain how the algorithm in \cite{cahnkjuchukova2016linking} is used to compute the linking numbers that appear in both examples.  The input for the algorithm is a labeled link diagram, which is recorded by means of several lists analogous to the Gauss code.  One component of the link diagram is the knot $\alpha$. In order to simplify the combinatorics, we only include  in our diagram two of $\{\beta, \omega_1, ... \omega_{2g-2}\}$, or one of these curves together with its push-off in $V$, at any given time.  Call these two curves $g$ and $h$.   Because $\beta$ is a mod 3 characteristic knot, any curve in $V-\beta$ lifts to three closed loops~\cite{CS1984linking}.  Thus for each pair of curves in $\mathcal{B}_V$, we compute nine linking numbers of their lifts, organized in a symmetric $3\times 3$ matrix. 

The following set-up allows us to compute the intersection number of any lift of $h$ with a 2-chain whose boundary is any given lift of $g$.  For the details on how this 2-chain is constructed see~\cite{cahnkjuchukova2016linking}.

{\it(1)} The arcs of $\alpha$ in the diagram $\alpha\cup g$ are labeled $0,1, \dots, m-1$, where $m$ is the number of self-crossings of $\alpha$ plus the number of crossings of $\alpha$ under $g$.  Each arc of $\alpha$ is colored 1,2 or 3, according to the given homomorphism $\rho:\pi_1(S^3-\alpha)\twoheadrightarrow D_3$.

{\it(2)} The arcs of $g$ in the diagram $\alpha\cup g $ are labeled $0,1,\dots, n-1$, where $n$ is the number of self-crossings of $g$ plus the number of crossings of $g$ under $\alpha$.

{\it(3)} Now we add $h$ to the above numbered diagram $\alpha\cup g$ without changing the numbering of any existing arcs.   The arcs of $h$ are labeled $0,1,\dots, o-1$, where $o$ is the number of crossings of $h$ under $\alpha$ plus the number of crossings of $h$ under $g$.  In this article, $h$ never has self-crossings. 

\subsection{Lists.}

We provide the input used in the computation of our Example~1, the knot $6_1$. The first four lists needed are associated to the knot $\alpha$.  The remaining six lists are associated to the two curves $g$ and $h$ described above.  The first list, $f=(f(i))_i$, records each number $f(i)$ assigned to the over-arc which meets the head of arc $i$ of $\alpha$.  The second, $\epsilon=(\epsilon(i))_i$, denotes the local writhe number at the head of arc $i$.  Next, $t=(t(i))_i$ denotes the {\it type} of crossing at the head of arc $i$; that is, we set $t(i)=p$ if the over-arc at the head of arc $i$ is an arc of $g$, and we set $t(i)=k$ if the over-arc at the head of arc $i$ is another arc of $\alpha$.  Recall that the $i^{th}$ arc of $\alpha$ may be a union of smaller arcs, separated by over-crossings by arcs of $h$. Due to our numbering system, the over-crossing at the end of an arc of $\alpha$ will never be an arc of $h$.  The fourth list, $c=(c(i))_i$, enumerates the colors on the consecutive arcs of $\alpha$.

Numbering, signs, crossing types, and colors for $\alpha$:
$$f=(1,8,0,7,10,5,3,2,4,6,6,4)$$
$$\epsilon=(-,+,-,-,-,-,+,+,+,-,+,-)$$
$$t=(p,k,k,p,k,p,p,k,p,k,p,k)$$
$$c=(1,1,3,2,2,1,1,1,2,2,3,3)$$

The remaining lists are the over-crossing numbers, signs, and crossing types for the other two components, $g$ and $h$, of the link diagram.

Numbering, signs, and crossing types for $\beta$:
$$(0,8,2,6,6,10,4,0)$$
$$(-,+,-,+,-,+,-,+)$$
$$(k,k,k,k,k,k,k,k)$$

Numbering, signs, and crossing types for $\beta_r$:
$$(0,1,8,2,3,6,4,6,10,6,4,0)$$
$$(-,-,+,-,-,+,+,-,+,+,-,+)$$
$$(k,p,k,k,p,k,p,k,k,p,k,k)$$\\

	{\bf Acknowledgment.} We would like to thank Julius Shaneson for helpful discussions. The anonymous referee provided valuable feedback. Parts of this work were completed at the Max Planck Institute for Mathematics. We are grateful to MPIM for its support and hospitality.  \\

Patricia Cahn\\
Smith College\\
{\it pcahn@smith.edu}\\

Alexandra Kjuchukova\\
Max Planck Institute for Mathematics -- Bonn\\
{\it sashka@mpim-bonn.mpg.de}

\bibliographystyle{amsplain}
\bibliography{BrCovBib}

\end{document}